\documentclass[12pt]{article}
\usepackage{amssymb}

\newtheorem{definition}{Definition}[section]
\newtheorem{theorem}[definition]{Theorem}
\newtheorem{lemma}[definition]{Lemma}
\newtheorem{corollary}[definition]{Corollary}
\newtheorem{remark}[definition]{Remark}

\newtheorem{note}[definition]{Note}

\def\K{\mathbb K}

\def\Z{\mathbb Z}
\def\K{\mathbb K}

\begin{document}

\title{ \bf The quantum
algebra $U_q(\mathfrak{sl}_2)$ 
and \\its equitable presentation\footnote{
{\bf Keywords}.  Quantum group, quantum algebra, Leonard pair,
tridiagonal pair.
 \hfil\break
\noindent {\bf 2000 Mathematics Subject Classification}. 
Primary: 17B37. Secondary: 20G42, 16W35, 81R50. 
}}
\author{Tatsuro Ito,
Paul Terwilliger and Chih-wen Weng
}
\date{}

\maketitle
\begin{abstract} 
We show that the quantum 
algebra $U_q(\mathfrak{sl}_2)$
has a presentation with generators $x^{\pm 1},y,z$
and relations $x x^{-1} = x^{-1} x=1$,
\begin{eqnarray*}
\frac{qxy-q^{-1}yx}{q-q^{-1}}=1,
\qquad 
\frac{qyz-q^{-1}zy}{q-q^{-1}}=1,
\qquad 
\frac{qzx-q^{-1}xz}{q-q^{-1}}=1.
\end{eqnarray*}
We call this the
{\it equitable} presentation. 
We show that $y$ (resp. $z$) is not invertible
in 
$U_q(\mathfrak{sl}_2)$ by
displaying an infinite dimensional
$U_q(\mathfrak{sl}_2)$-module that contains
a nonzero null vector for
$y$ (resp. $z$). 
We consider finite dimensional
$U_q(\mathfrak{sl}_2)$-modules under
the assumption that
$q$ is not a root of 1 and 
$\mbox{char}(\K) \not=2$, where
$\K$ is the underlying field.
We  
show that $y$ and $z$ are invertible on
each finite dimensional 
$U_q(\mathfrak{sl}_2)$-module.
We display a linear operator $\Omega$ that acts
on finite dimensional 
$U_q(\mathfrak{sl}_2)$-modules, and satisfies
\begin{eqnarray*}
\Omega^{-1}x\, \Omega = y,
\qquad 
\Omega^{-1}y\,  \Omega = z,
\qquad 
\Omega^{-1}z\, \Omega = x
\end{eqnarray*}
on these modules. We define $\Omega$ using the $q$-exponential
function.
\end{abstract}

\section{The  algebra 
$U_q(\mathfrak{sl}_2)$
}

Let $\K$ denote a field
and let $q$ denote
a nonzero scalar in $\K$ such that $q^2\not=1$.
For an integer $n$ we define
\begin{eqnarray*}
\lbrack n \rbrack = \frac {q^n-q^{-n}}{q-q^{-1}}
\end{eqnarray*}
and for $n\geq 0$ we define
\begin{eqnarray*}
\lbrack n \rbrack^{!} = 
\lbrack n \rbrack  \lbrack n-1 \rbrack
\cdots \lbrack 2 \rbrack \lbrack 1 \rbrack.
\end{eqnarray*}
We interpret $\lbrack 0 \rbrack^!=1$. 
We now recall
the quantum algebra $U_q(\mathfrak{sl}_2)$.

\begin{definition} 
\label{def:uq}
\rm
We let $U_q(\mathfrak{sl}_2)$
denote the unital associative $\K$-algebra 
with
generators $k^{\pm 1}$,$e$, $f$
and the following relations:
\begin{eqnarray*}
kk^{-1} &=& 
k^{-1}k =  1,
\label{eq:buq1}
\\
ke &=& q^2ek,
\label{eq:buq2}
\\
kf&=&q^{-2}fk,
\label{eq:buq3}
\\
ef-fe &=& \frac{k-k^{-1}}{q-q^{-1}}.
\label{eq:buq4}
\end{eqnarray*}
\end{definition}

\noindent 
We call $k^{\pm 1}, e,f$
the {\it Chevalley generators} for
$U_q(\mathfrak{sl}_2)$.

\medskip
\noindent 
We refer the reader to
 \cite{Jantzen},
\cite{Kassel}
 for
background information on 
$U_q(\mathfrak{sl}_2)$. We will generally follow the
notational conventions of \cite{Jantzen}.

%

\section{The equitable presentation for
$U_q(\mathfrak{sl}_2)$
}

In
the presentation for 
$U_q(\mathfrak{sl}_2)$ given
in 
Definition
\ref{def:uq}
the generators $k^{\pm 1}$
and the generators $e,f$ play a very different
role. 
We now introduce a presentation for 
$U_q(\mathfrak{sl}_2)$ whose
generators
are on a more equal footing.  

\begin{theorem}
\label{thm:uq2}
The algebra
$U_q(\mathfrak{sl}_2)$ 
 is isomorphic to
the unital associative $\K$-algebra 
with
generators 
$x^{\pm 1}$, $y$, $z$
and the following relations:
\begin{eqnarray}
xx^{-1} = 
x^{-1}x &=&  1,
\label{eq:2buq1}
\\
\frac{qxy-q^{-1}yx}{q-q^{-1}}&=&1,
\label{eq:2buq2}
\\
\frac{qyz-q^{-1}zy}{q-q^{-1}}&=&1,
\label{eq:2buq3}
\\
\frac{qzx-q^{-1}xz}{q-q^{-1}}&=&1.
\label{eq:2buq4}
\end{eqnarray}
An isomorphism with the presentation in Definition
\ref{def:uq} is given by:
\begin{eqnarray*}
\label{eq:iso1}
x^{{\pm}1} &\rightarrow & k^{{\pm}1},\\
\label{eq:iso2}
y &\rightarrow & k^{-1}+f(q-q^{-1}), \\
\label{eq:iso3}
z &\rightarrow & k^{-1}-k^{-1}eq(q-q^{-1}).
\end{eqnarray*}
The inverse of this isomorphism is given by:
\begin{eqnarray*}
\label{eq:iso1inv}
k^{{\pm}1} &\rightarrow & x^{{\pm}1},\\
\label{eq:iso2inv}
f &\rightarrow & (y-x^{-1})(q-q^{-1})^{-1}, \\
\label{eq:iso3inv}
e &\rightarrow & (1-xz)q^{-1}(q-q^{-1})^{-1}.
\end{eqnarray*}
\end{theorem}
\noindent {\it Proof:} One readily checks that
each map is a homomorphism of $\K$-algebras and that
the maps are inverses.
It follows that each map is an isomorphism of $\K$-algebras.
\hfill $\Box $ \\

\noindent
The generators 
 $x^{\pm 1}$, $y$, $z$ from
Theorem
\ref{thm:uq2} 
 are 
on an equal footing, more or less.
In view
of this we make a definition.

\begin{definition}
\rm
By the {\it equitable presentation} for
$U_q(\mathfrak{sl}_2)$ we mean
the presentation given in
 Theorem
\ref{thm:uq2}.
We call 
$x^{\pm 1}, y,z$ the
{\it equitable generators}.
\end{definition}

\medskip
\noindent 
We remark that the isomorphism given in Theorem
\ref{thm:uq2} is not unique. This is a consequence of
the following lemma.  

\begin{lemma}
\label{lem:aut}
For an integer $i$ and nonzero 
$\alpha \in \K$ there
exists a
$\K$-algebra automorphism 
of $U_q(\mathfrak{sl}_2)$ that satisfies
\begin{eqnarray*}
k^{\pm 1} \rightarrow k^{\pm 1},
\qquad \qquad 
e \rightarrow \alpha e k^i,
\qquad \qquad 
f \rightarrow \alpha^{-1}k^{-i}f.
\end{eqnarray*}
\end{lemma}
\noindent {\it Proof:} 
Routine.
\hfill $\Box $ \\

\section{The elements $y$ and $z$ are not invertible in 
$U_q(\mathfrak{sl}_2)$}

\noindent In this section we show that
the equitable generators
$y$ and $z$ are not invertible in
$U_q(\mathfrak{sl}_2)$. 
In order to show that $y$ (resp. $z$)
is not invertible
in 
$U_q(\mathfrak{sl}_2)$
we 
display an infinite dimensional 
$U_q(\mathfrak{sl}_2)$-module
that contains a nonzero null vector for $y$ (resp. $z$).

\begin{lemma}
\label{lem:mody} There exists a
$U_q(\mathfrak{sl}_2)$-module $\Gamma_y$ with the following
property: $\Gamma_y$ has a basis 
\begin{eqnarray*}
u_{ij} \qquad i,j\in \Z, \quad j\geq 0
\end{eqnarray*}
such that
\begin{eqnarray*}
x u_{ij} &=& u_{i+1,j},
\\
x^{-1} u_{ij} &=& u_{i-1,j},
\\
y u_{ij} &=& 
q^{2i-j}(q^j-q^{-j})u_{i,j-1}-
q^i(q^i-q^{-i})u_{i-1,j},
\\
z u_{ij} &=& q^{-2i}u_{i,j+1}+q^{-i}(q^i-q^{-i})u_{i-1,j}
\end{eqnarray*}
for all $i,j \in \Z$ with $j\geq 0$. In the above equations
$u_{r,-1}:=0$ for $r \in \Z$.
\end{lemma}
\noindent {\it Proof:} 
One routinely verifies that the given actions of $x^{\pm 1}, y,z$
satisfy the  relations
(\ref{eq:2buq1})--(\ref{eq:2buq4}).
\hfill $\Box $ \\

\begin{lemma} The following {\rm (i)--(iii)} hold.
\begin{enumerate}
\item[{\rm (i)}] $yu_{00}=0$, where the vector $u_{00}$ is from
Lemma
\ref{lem:mody}.
\item[{\rm (ii)}]
$y$ is not invertible on $\Gamma_y$, where $\Gamma_y$ is
the 
$U_q(\mathfrak{sl}_2)$-module from
Lemma
\ref{lem:mody}.
\item[{\rm (iii)}]
 $y$ is not invertible in 
$U_q(\mathfrak{sl}_2)$.
\end{enumerate}
\end{lemma}
\noindent {\it Proof:}  Immediate.
\hfill $\Box $ \\

\begin{remark}
\rm Referring to Lemma
\ref{lem:mody}, we have $u_{ij}= x^iz^ju_{00}$ for
$i, j\in \Z$, $j\geq 0$. 
\end{remark}

\begin{lemma}
\label{lem:mod1} There exists a
$U_q(\mathfrak{sl}_2)$-module $\Gamma_z$ with the following
property: $\Gamma_z$ has a basis 
\begin{eqnarray*}
v_{ij} \qquad i,j\in \Z, \quad j\geq 0
\end{eqnarray*}
such that
\begin{eqnarray*}
x v_{ij} &=& v_{i+1,j},
\\
x^{-1} v_{ij} &=& v_{i-1,j},
\\
y v_{ij} &=& q^{2i}v_{i,j+1}-q^i(q^i-q^{-i})v_{i-1,j},
\\
z v_{ij} &=& q^{-i}(q^i-q^{-i})v_{i-1,j}-q^{j-2i}(q^j-q^{-j})v_{i,j-1}
\end{eqnarray*}
for all $i,j \in \Z$ with $j\geq 0$. In the above equations
$v_{r,-1}:=0$ for $r \in \Z$.
\end{lemma}
\noindent {\it Proof:} 
One routinely verifies that the given actions of $x^{\pm 1}, y,z$
satisfy the  relations
(\ref{eq:2buq1})--(\ref{eq:2buq4}).
\hfill $\Box $ \\

\begin{lemma} The following {\rm (i)--(iii)} hold.
\begin{enumerate}
\item[{\rm (i)}]
$zv_{00}=0$, where the vector $v_{00}$ is from
Lemma
\ref{lem:mod1}.
\item[{\rm (ii)}]
$z$ is not invertible on $\Gamma_z$,
 where $\Gamma_z$ is
the 
$U_q(\mathfrak{sl}_2)$-module from
Lemma
\ref{lem:mod1}.
\item[{\rm (iii)}]
 $z$ is not invertible in 
$U_q(\mathfrak{sl}_2)$.
\end{enumerate}
\end{lemma}
\noindent {\it Proof:}  Immediate.
\hfill $\Box $ \\

\begin{remark}
\rm Referring to Lemma
\ref{lem:mod1}, we have $v_{ij}= x^iy^jv_{00}$ for
$i, j\in \Z$, $j\geq 0$. 
\end{remark}

\section{Finite dimensional 
$U_q(\mathfrak{sl}_2)$-modules}

\noindent 
From now on we restrict
our attention to finite
dimensional 
$U_q(\mathfrak{sl}_2)$-modules.
In order to simplify things we 
make the following assumption.
\begin{eqnarray*}
\mbox{
For the rest of this paper,
we assume $q$ is not a root of 1, and that
$\mbox{char}(\K) \not=2$.
}
\end{eqnarray*}
In this section we
show that the equitable generators $y$ and $z$ 
are invertible on each finite dimensional
$U_q(\mathfrak{sl}_2)$-module.

\medskip
\noindent 
We begin with some general comments.
By \cite[Theorems 2.3, 2.9]{Jantzen} each finite dimensional
$U_q(\mathfrak{sl}_2)$-module $M$
is semi-simple; this means that $M$
is a direct sum of simple
$U_q(\mathfrak{sl}_2)$-modules.
The finite dimensional simple 
$U_q(\mathfrak{sl}_2)$-modules
are described as follows.

\begin{lemma} \cite[Theorem 2.6]{Jantzen}
\label{lem:lneps}
There exists a family of finite dimensional
simple 
$U_q(\mathfrak{sl}_2)$-modules
\begin{eqnarray}
L(n,\varepsilon) \qquad \qquad 
\varepsilon \in \lbrace 1,-1\rbrace,
\qquad \qquad 
n=0,1,2,\ldots 
\label{eq:smods}
\end{eqnarray}
with the following properties: 
$L(n,\varepsilon)$
has a basis $v_0, v_1, \ldots, v_n$
such that $kv_i =\varepsilon q^{n-2i}v_i$
for $0 \leq i \leq n$,
$fv_i = \lbrack i+1 \rbrack v_{i+1}$ for $0 \leq i \leq n-1$,
$fv_n=0$,
$ev_i = \varepsilon\lbrack n-i+1 \rbrack v_{i-1}$ for 
$1 \leq i \leq n$,
$ev_0=0$.
Every finite dimensional simple 
$U_q(\mathfrak{sl}_2)$-module is isomorphic to
exactly one of
the modules (\ref{eq:smods}).
\end{lemma}

\noindent 
The equitable generators act on
the modules
$L(n,\varepsilon)$ as follows.

\begin{lemma}
\label{lem:altact}
For an integer $n\geq 0$ and for
$\varepsilon \in \lbrace 1, -1\rbrace $,
the 
$U_q(\mathfrak{sl}_2)$-module $L(n,\varepsilon)$
has a basis $u_0, u_1, \ldots, u_n$ such that
\begin{eqnarray}
\varepsilon x u_i &=& q^{n-2i}u_i \qquad \qquad
(0 \leq i \leq n),
\label{eq:xact}
\\
(\varepsilon y -q^{2i-n}) u_i &=& (q^{-n}-q^{2i+2-n})u_{i+1}
\qquad
(0 \leq i \leq n-1),
\qquad (\varepsilon y -q^n)u_n=0,
\label{eq:yact}
\\
(\varepsilon z -q^{2i-n}) u_i &=& (q^n-q^{2i-2-n})u_{i-1}
\qquad
(1 \leq i \leq n),
\qquad (\varepsilon z -q^{-n})u_0=0.
\label{eq:zact}
\end{eqnarray}
\end{lemma}
\noindent {\it Proof:} 
For the purpose of this proof we identify the copy
of
$U_q(\mathfrak{sl}_2)$
given in Definition
\ref{def:uq} with the copy
given in Theorem
\ref{thm:uq2}, via the isomorphism
in Theorem
\ref{thm:uq2}.
Let the basis $v_0,v_1, \ldots, v_n$
for $L(n,\varepsilon)$ be as in
Lemma
\ref{lem:lneps}.
Define $u_i = \gamma_i v_i$
for $0 \leq i \leq n$, where
$\gamma_0= 1$ and $\gamma_i = -\varepsilon q^{n-i}\gamma_{i-1}
$ for $1 \leq i \leq n$. Using $x=k$,
$y=k^{-1}+f(q-q^{-1})$, and $z=k^{-1}-k^{-1}eq(q-q^{-1})$,
together with the data in Lemma
\ref{lem:lneps}, we routinely verify
(\ref{eq:xact})--(\ref{eq:zact}).
\hfill $\Box $ \\

\begin{note} \rm
The basis $u_0, u_1, \ldots, u_n$ in Lemma
\ref{lem:altact} is normalized so that
$y u = \varepsilon q^{-n}u$
and 
$z u = \varepsilon q^{n}u$
for $u=\sum_{i=0}^n u_i$.
\end{note}

\begin{corollary}
\label{cor:inv}
For an integer $n\geq 0$ and for
 $\varepsilon \in \lbrace 1,-1\rbrace $,
the following {\rm (i), (ii)} hold on the
$U_q(\mathfrak{sl}_2)$-module $L(n,\varepsilon)$.
\begin{enumerate}
\item[{\rm (i)}]
Each of $x,y,z$ is semi-simple with eigenvalues
$\varepsilon q^n, 
\varepsilon q^{n-2}, 
\ldots,
\varepsilon q^{-n}$.
\item[{\rm (ii)}]
Each of $x,y,z$ is invertible.
\end{enumerate}
\end{corollary}
\noindent {\it Proof:} 
For $x$ this is clear from
(\ref{eq:xact}). We now verify our assertions for $y$.
With respect to the basis 
$u_0, u_1, \ldots, u_n$ for $L(n,\varepsilon)$
given in Lemma
\ref{lem:altact}, 
by (\ref{eq:yact}) the matrix representing
$y$ is lower triangular with $(i,i)$ entry 
$\varepsilon q^{2i-n}$ for $0 \leq i \leq n$.
Therefore the action of $y$ on $L(n,\varepsilon)$
has eigenvalues 
$\varepsilon q^n, 
\varepsilon q^{n-2}, 
\ldots,
\varepsilon q^{-n}$.
These eigenvalues are mutually distinct so this action
is semi-simple. These eigenvalues are nonzero so
this action is invertible.
We have now verified our assertions for $y$. Our
assertions for $z$ are similarly verified.
\hfill $\Box $ \\

\noindent By Corollary
\ref{cor:inv} and since each finite dimensional 
$U_q(\mathfrak{sl}_2)$-module 
is semi-simple we obtain the following result.

\begin{corollary}
\label{cor:yzinv}
On each finite dimensional 
$U_q(\mathfrak{sl}_2)$-module 
the actions of $y$ and $z$ are invertible.
\end{corollary}

\noindent Motivated by Corollary
\ref{cor:yzinv}
we make the following definition.

\begin{definition}
\label{def:invop}
\rm
We let $y^{-1}$ (resp. $z^{-1}$) denote the
linear operator that acts on each finite
dimensional 
$U_q(\mathfrak{sl}_2)$-module 
as the inverse of $y$ (resp. $z$).
\end{definition}

\section{The elements $n_x, n_y, n_z$}

In this section we define some elements 
$n_x, n_y, n_z$ of
$U_q(\mathfrak{sl}_2)$ and show
that these are nilpotent on each finite
dimensional 
$U_q(\mathfrak{sl}_2)$-module.
We then recall the $q$-exponential
function ${\rm{exp}}_q$ and derive a number
of equations involving 
${\rm{exp}}_q(n_x),
{\rm{exp}}_q(n_y),
{\rm{exp}}_q(n_z)$.
These equations will show that on finite dimensional
$U_q(\mathfrak{sl}_2)$-modules the operators
$y^{-1}, z^{-1}$ from Definition
\ref{def:invop} satisfy
\begin{eqnarray}
y^{-1}&=& {\rm{exp}}_q(n_z) \,x \,{\rm{exp}}_q(n_z)^{-1},
\label{eq:yinvint}
\\
z^{-1}&=& {\rm{exp}}_q(n_y)^{-1}\,x \,{\rm{exp}}_q(n_y).
\label{eq:zinvint}
\end{eqnarray}

\noindent We begin with an observation. 

\begin{lemma}
\label{lem:relrewrite}
The equitable generators $x,y,z$ of 
$U_q(\mathfrak{sl}_2)$ 
satisfy
\begin{eqnarray*}
 q(1-yz) &=& 
 q^{-1}(1-zy),
\label{eq:xalt}
\\
q(1-zx) &=& 
q^{-1}(1-xz),
\label{eq:yalt}
\\
q(1-xy) &=& 
q^{-1}(1-yx).
\label{eq:zalt}
\end{eqnarray*}
\end{lemma}
\noindent {\it Proof:} These equations
are reformulations
of
(\ref{eq:2buq2})--(\ref{eq:2buq4}).
\hfill $\Box $ \\

\begin{definition}
\label{def:n}
\rm
We let  $n_x, n_y, n_z$ denote the following elements
in
$U_q(\mathfrak{sl}_2)$:
\begin{eqnarray}
n_x &=& \frac{q(1-yz)}{q-q^{-1}} = 
 \frac{q^{-1}(1-zy)}{q-q^{-1}},
\label{eq:nx}
\\
n_y &=& \frac{q(1-zx)}{q-q^{-1}} = 
 \frac{q^{-1}(1-xz)}{q-q^{-1}},
\label{eq:ny}
\\
n_z &=& \frac{q(1-xy)}{q-q^{-1}} = 
 \frac{q^{-1}(1-yx)}{q-q^{-1}}.
\label{eq:nz}
\end{eqnarray}
\end{definition}

\begin{note}
\rm
Under the isomorphism given in Theorem
\ref{thm:uq2} the preimage of $n_y$ (resp. $n_z$)
is $e$ (resp. $-qkf$).
\end{note}

\noindent We recall some notation.
Let $V$ denote a finite dimensional vector space
over $\K$. A linear transformation $T:V\rightarrow V$
is called {\it nilpotent} whenever there exists a positive
integer $r$ such that $T^rV=0$.

\medskip
\noindent We are going to show that each of $n_x, n_y, n_z$
is nilpotent on all finite dimensional 
$U_q(\mathfrak{sl}_2)$-modules. We will show this
using the following lemma.

\begin{lemma}
\label{lem:nrel}
The following relations hold in
$U_q(\mathfrak{sl}_2)$:
\begin{eqnarray}
xn_y &=& q^2 n_y x, 
\qquad \qquad 
xn_z = q^{-2} n_z x,
\label{eq:xn}
\\
yn_z &=& q^2 n_z y, 
\qquad \qquad 
yn_x = q^{-2} n_x y,
\label{eq:yn}
\\
zn_x &=& q^2 n_x z, 
\qquad \qquad 
zn_y = q^{-2} n_y z.
\label{eq:zn}
\end{eqnarray}
\end{lemma}
\noindent {\it Proof:} 
In order to verify these equations, 
eliminate $n_x, n_y, n_z$ using
Definition \ref{def:n} and simplify
the result.
\hfill $\Box $ \\

\begin{lemma}
\label{lem:allnil}
Each of $n_x, n_y, n_z$ is nilpotent on finite
dimensional 
$U_q(\mathfrak{sl}_2)$-modules.
\end{lemma}
\noindent {\it Proof:}
We prove the result for $n_x$; the proof for
$n_y$ and $n_z$ is similar.
Since each finite dimensional 
$U_q(\mathfrak{sl}_2)$-module is semi-simple
and in view of Lemma
\ref{lem:lneps}, it suffices to show that
$n_x$ is nilpotent on each module
$L(n,\varepsilon)$. By Corollary
\ref{cor:inv}(i),  $L(n,\varepsilon)$ has a basis
$w_0, w_1, \ldots, w_n$ such that
$yw_i = \varepsilon q^{n-2i}w_i$
for $0 \leq i \leq n$. Using the equation
on the right in
(\ref{eq:yn}) we routinely find
that $n_x w_i$ is a scalar multiple of
$w_{i+1}$ for $0 \leq i \leq n-1$
and $n_x w_n=0$. 
This shows that
$n_x$ is nilpotent on
$L(n,\varepsilon)$ and the result follows.
\hfill $\Box $ \\

\noindent We now recall the $q$-exponential function.

\begin{definition} \cite[p.~204]{tanis}
\label{def:tnil}
\rm
Let $T$ denote a linear operator that acts on
finite dimensional 
$U_q(\mathfrak{sl}_2)$-modules in a nilpotent
fashion.
We define
\begin{eqnarray}
{\rm{exp}}_q(T) &=& \sum_{i=0}^\infty 
\frac{q^{i(i-1)/2}}{\lbrack i \rbrack^!}\, T^i.
\label{eq:qexp}
\end{eqnarray}
We view ${\rm{exp}}_q(T)$ as a linear operator
that acts on finite dimensional 
$U_q(\mathfrak{sl}_2)$-modules.
\end{definition}

\noindent The following result is well known and
easily verified.

\begin{lemma} \cite[p.~204]{tanis}
\label{lem:expinv}
Let $T$ denote a linear operator that acts on
finite dimensional 
$U_q(\mathfrak{sl}_2)$-modules in a nilpotent
fashion.
Then on each of these modules
${\rm{exp}}_q(T)$ is
invertible; the inverse is
\begin{eqnarray*}
{\rm{exp}}_{q^{-1}}(-T) =\sum_{i=0}^\infty 
\frac{(-1)^i q^{-i(i-1)/2}}{\lbrack i \rbrack^!}\,T^i.
\end{eqnarray*}
\end{lemma}

\noindent The following lemma is essentially a special case of
\cite[Eq. (4.8)]{ros}; we include a proof
for the sake of completeness. 

\begin{lemma} 
\cite[Eq. (4.8)]{ros}
\label{lem:expform2}
The following {\rm (i)--(iii)} hold on each finite
dimensional 
$U_q(\mathfrak{sl}_2)$-module.
\begin{enumerate}
\item[{\rm (i)}]
${\rm{exp}}_q(n_y)^{-1}\,x\, {\rm{exp}}_q(n_y) = z^{-1}$,
\item[{\rm (ii)}]
${\rm{exp}}_q(n_z)^{-1}\,y\, {\rm{exp}}_q(n_z) = x^{-1}$,
\item[{\rm (iii)}]
${\rm{exp}}_q(n_x)^{-1}\,z\, {\rm{exp}}_q(n_x) = y^{-1}$.
\end{enumerate}
\end{lemma}
\noindent {\it Proof:} 
(i) We show
\begin{eqnarray}
x\,{\rm exp}_q(n_y) \,z = {\rm exp}_q(n_y).
\label{eq:nzx}
\end{eqnarray}
The left side of
(\ref{eq:nzx}) is equal to
$x\,{\rm exp}_q(n_y)\, x^{-1}xz$.
Observe
$x\,{\rm exp}_q(n_y) \,x^{-1} = {\rm exp}_q(xn_yx^{-1})$
by
(\ref{eq:qexp}) and 
$xn_yx^{-1} = q^2n_y$ by 
(\ref{eq:xn}). Also
$xz=1-q(q-q^{-1})n_y$ by
(\ref{eq:ny}). Using these comments
and 
(\ref{eq:qexp}) we routinely find
that the left side of 
(\ref{eq:nzx}) is equal to
the right side of 
(\ref{eq:nzx}). The result follows.
\\
\noindent (ii), (iii) Similar to the proof of (i) above.
\hfill $\Box $ \\

\noindent For convenience we display a second ``version''
of Lemma
\ref{lem:expform2}.

\begin{lemma}
\label{lem:expform}
The following {\rm (i)--(iii)} hold on each finite
dimensional 
$U_q(\mathfrak{sl}_2)$-module.
\begin{enumerate}
\item[{\rm (i)}]
${\rm{exp}}_q(n_z)\,x\, {\rm{exp}}_q(n_z)^{-1} = y^{-1}$,
\item[{\rm (ii)}]
${\rm{exp}}_q(n_x)\,y\, {\rm{exp}}_q(n_x)^{-1} = z^{-1}$,
\item[{\rm (iii)}]
${\rm{exp}}_q(n_y)\,z\, {\rm{exp}}_q(n_y)^{-1} = x^{-1}$.
\end{enumerate}
\end{lemma}
\noindent {\it Proof:} 
For each of the equations in
Lemma
\ref{lem:expform2} take the inverse of each side
and simplify the result.
\hfill $\Box $ \\

\noindent We note that
the equations
(\ref{eq:yinvint}), 
(\ref{eq:zinvint})  are just
Lemma
\ref{lem:expform}(i)
and Lemma
\ref{lem:expform2}(i),
 respectively.

\section{Some formulae involving the $q$-exponential function}

\noindent In the next section we will display a linear
operator $\Omega$ that acts on finite dimensional
$U_q(\mathfrak{sl}_2)$-modules, and satisfies
$\Omega^{-1}x\, \Omega = y$,
$\Omega^{-1}y\,  \Omega = z$,
$\Omega^{-1}z\, \Omega = x$ on these modules.
In order to prove that $\Omega $ has the desired
properties we will first establish
a few identities. These identities are given
in this section.

\begin{lemma}
\label{lem:xyx}
The following {\rm (i)--(iii)} hold on each finite
dimensional 
$U_q(\mathfrak{sl}_2)$-module.
\begin{enumerate}
\item[{\rm (i)}]
${\rm{exp}}_q(n_z)^{-1}\,x\, {\rm{exp}}_q(n_z) = xyx$,
\item[{\rm (ii)}]
${\rm{exp}}_q(n_x)^{-1}\,y\, {\rm{exp}}_q(n_x) = yzy$,
\item[{\rm (iii)}]
${\rm{exp}}_q(n_y)^{-1}\,z\, {\rm{exp}}_q(n_y) = zxz$.
\end{enumerate}
\end{lemma}
\noindent {\it Proof:}
(i) The element $xy$ commutes with $n_z$ by
(\ref{eq:nz}) so $xy$ commutes with 
${\rm exp}_q(n_z)$ in view of
(\ref{eq:qexp}).
Therefore
${\rm exp}_q(n_z)^{-1} \,xy\, {\rm exp}_q(n_z)=xy$.
By Lemma 
\ref{lem:expform2}(ii) we have
$y \, {\rm exp}_q(n_z)={\rm exp}_q(n_z)\, x^{-1}$.
Combining these last two equations we routinely
obtain the result.
\\
\noindent (ii), (iii) Similar to the proof of (i) above.
\hfill $\Box $ \\

\begin{lemma}
\label{lem:xyx2}
The following {\rm (i)--(iii)} hold on each finite
dimensional 
$U_q(\mathfrak{sl}_2)$-module.
\begin{enumerate}
\item[{\rm (i)}]
${\rm{exp}}_q(n_y)\,x\, {\rm{exp}}_q(n_y)^{-1} = xzx$,
\item[{\rm (ii)}]
${\rm{exp}}_q(n_z)\,y\, {\rm{exp}}_q(n_z)^{-1} = yxy$,
\item[{\rm (iii)}]
${\rm{exp}}_q(n_x)\,z\, {\rm{exp}}_q(n_x)^{-1} = zyz$.
\end{enumerate}
\end{lemma}
\noindent {\it Proof:}
(i) By Lemma
\ref{lem:xyx}(iii) we have
${\rm{exp}}_q(n_y)\,zxz\, {\rm{exp}}_q(n_y)^{-1} = z$.
In this equation we eliminate 
${\rm{exp}}_q(n_y)\,z$ and
$z\, {\rm{exp}}_q(n_y)^{-1}$
using Lemma
\ref{lem:expform}(iii). The result follows.
\\
\noindent (ii), (iii) Similar to the proof of
(i) above.
\hfill $\Box $ \\

\begin{lemma}
\label{lem:xyy}
The following {\rm (i)--(iii)} hold on each finite
dimensional 
$U_q(\mathfrak{sl}_2)$-module.
\begin{enumerate}
\item[{\rm (i)}]
${\rm{exp}}_q(n_x)^{-1}\,x\, {\rm{exp}}_q(n_x) = x+y-y^{-1}$,
\item[{\rm (ii)}]
${\rm{exp}}_q(n_y)^{-1}\,y\, {\rm{exp}}_q(n_y) = y+z-z^{-1}$,
\item[{\rm (iii)}]
${\rm{exp}}_q(n_z)^{-1}\,z\, {\rm{exp}}_q(n_z) = z+x-x^{-1}$.
\end{enumerate}
\end{lemma}
\noindent {\it Proof:}
(i) Using 
(\ref{eq:2buq2}),
(\ref{eq:2buq4})
and 
 (\ref{eq:nx})
we obtain 
$xn_x-n_x x = y-z$. By this and a routine induction
using
(\ref{eq:yn}), 
(\ref{eq:zn})
we find
\begin{eqnarray}
\label{eq:nxind}
x n_x^i - n_x^i x= q^{1-i} \lbrack i \rbrack (n_x^{i-1} y - z n_x^{i-1})
\end{eqnarray}
for each integer $i\geq 0$. Using 
(\ref{eq:qexp}) and
(\ref{eq:nxind})
we obtain
\begin{eqnarray}
\label{eq:nxindexp}
x\,
{\rm exp}_q(n_x)-{\rm exp}_q(n_x)\, x = {\rm exp}_q(n_x)\,y
-z \,{\rm exp}_q(n_x).
\end{eqnarray}
In line 
(\ref{eq:nxindexp}) we multiply each term on
the left by 
${\rm exp}_q(n_x)^{-1}$ and evaluate the term containing
$z$ using Lemma
\ref{lem:expform2}(iii) to get the result.
\\
\noindent (ii), (iii) Similar to the proof of (i) above.
\hfill $\Box $ \\

\begin{lemma}
\label{lem:xyy2}
The following {\rm (i)--(iii)} hold on each finite
dimensional 
$U_q(\mathfrak{sl}_2)$-module.
\begin{enumerate}
\item[{\rm (i)}]
${\rm{exp}}_q(n_x)\,x\, {\rm{exp}}_q(n_x)^{-1} = x+z-z^{-1}$,
\item[{\rm (ii)}]
${\rm{exp}}_q(n_y)\,y\, {\rm{exp}}_q(n_y)^{-1} = y+x-x^{-1}$,
\item[{\rm (iii)}]
${\rm{exp}}_q(n_z)\,z\, {\rm{exp}}_q(n_z)^{-1} = z+y-y^{-1}$.
\end{enumerate}
\end{lemma}
\noindent {\it Proof:}
(i) By Lemma
\ref{lem:xyy}(i) we have
\begin{eqnarray}
\label{eq:vvv}
x ={\rm{exp}}_q(n_x) \,(x+y-y^{-1})\, {\rm{exp}}_q(n_x)^{-1}.
\end{eqnarray}
By Lemma
\ref{lem:expform}(ii)
we have
${\rm{exp}}_q(n_x)\,y\,{\rm{exp}}_q(n_x)^{-1}=z^{-1}$
and 
${\rm{exp}}_q(n_x)\,y^{-1}\,{\rm{exp}}_q(n_x)^{-1}=z$.
Evaluating
(\ref{eq:vvv}) using these comments
we obtain the result.
\\
\noindent (ii), (iii) Similar to the proof of
(i) above.
\hfill $\Box $ \\

\section{The operator $\Omega$}

\noindent In this section we display a linear
operator $\Omega$ that acts on finite dimensional
$U_q(\mathfrak{sl}_2)$-modules, and satisfies
$\Omega^{-1}x\, \Omega = y$,
$\Omega^{-1}y\,  \Omega = z$,
$\Omega^{-1}z\, \Omega = x$ on these modules.
In order to define $\Omega$
we first recall the notion of a {\it weight space}.

\begin{definition}
\label{def:ws}
\rm
Let $M$ denote a finite dimensional 
$U_q(\mathfrak{sl}_2)$-module.
For an integer  $\lambda $ and for $\varepsilon \in \lbrace 1,-1 \rbrace$
define
\begin{eqnarray*}
M(\varepsilon,\lambda) = \lbrace v \in M \; \vert \;xv=\varepsilon 
q^\lambda v \rbrace.
\end{eqnarray*}
We call $M(\varepsilon, \lambda)$ the 
 $(\varepsilon, \lambda)$-{\it weight space}
of $M$ with respect to $x$. 
By Corollary
\ref{cor:inv}(i) and since $M$ is semi-simple,
$M$ is the direct sum of its weight spaces with
respect to $x$.
\end{definition}

\begin{definition}
\label{def:psi}
\rm
We define a linear operator $\Psi$ that acts on
each finite dimensional 
$U_q(\mathfrak{sl}_2)$-module $M$. In order to
do this we give the action of $\Psi$ on each weight
space of $M$ with respect to $x$. For an integer
$\lambda $ and for $\varepsilon \in \lbrace 1,-1 \rbrace$,
$\Psi$ acts on the weight space $M(\varepsilon, \lambda)$
as $q^{-\lambda^2/2} I$ (if $\lambda $ is even)
and 
$q^{(1-\lambda^2)/2} I$ (if $\lambda $ is odd), where
$I$ denotes the identity map.
We observe that $\Psi$ is invertible on $M$.
\end{definition}

\begin{lemma}
\label{lem:psiprop}
For the operator $\Psi$ from Definition
\ref{def:psi} the following {\rm (i)--(iii)} hold on each
finite dimensional 
$U_q(\mathfrak{sl}_2)$-module.
\begin{enumerate}
\item[{\rm (i)}]
$\Psi^{-1}x \Psi=x$,
\item[{\rm (ii)}]
$\Psi^{-1}n_y \Psi=xn_y x$,
\item[{\rm (iii)}]
$\Psi^{-1}n_z \Psi=x^{-1}n_z x^{-1}$.
\end{enumerate}
\end{lemma}
\noindent {\it Proof:} 
Let $M$ denote a finite dimensional 
$U_q(\mathfrak{sl}_2)$-module.
For an integer $\lambda $ and for
 $\varepsilon \in \lbrace 1,-1 \rbrace$ we show
 that each of (i)--(iii) holds on $M(\varepsilon, \lambda)$.
\\
\noindent (i) On 
  $M(\varepsilon, \lambda)$ each of $\Psi, x$ acts as a scalar
  multiple of the identity.
\\
\noindent (ii) For notational convenience define
$s=0$ (if $\lambda $ is even) and $s=1$ (if $\lambda$ is odd).
For $v\in 
  M(\varepsilon, \lambda)$
we show $\Psi^{-1} n_y \Psi v = xn_yxv$.
Using the equation on the left in 
(\ref{eq:xn}) we find $n_y v \in M(\varepsilon, \lambda +2)$. Using
this we find
\begin{eqnarray*}
\Psi^{-1} n_y \Psi v &=& q^{(s-\lambda^2)/2} \Psi^{-1} n_y v
\\
&=&  q^{(s-\lambda^2)/2} 
 q^{((\lambda+2)^2-s)/2} 
n_y v
\\
&=&  q^{2 \lambda +2} 
n_y v
\end{eqnarray*}
and also
\begin{eqnarray*}
x n_y x v &=& \varepsilon q^{\lambda} x n_y v
\\ 
 &=& \varepsilon q^{\lambda} 
  \varepsilon q^{\lambda+2} 
 n_y v
\\
 &=& q^{2\lambda+2} 
 n_y v.
\end{eqnarray*}
Therefore 
$\Psi^{-1} n_y \Psi v = xn_yxv$.
We have now shown
$\Psi^{-1} n_y \Psi$ and $xn_yx$ coincide on
  $M(\varepsilon, \lambda)$.
\\
\noindent (iii) Similar to the proof of (ii) above.
\hfill $\Box $ \\

\begin{definition}
\label{def:omega}
\rm
We define
\begin{eqnarray}
\Omega = 
{\rm{exp}}_q(n_z)\, \Psi \,
{\rm{exp}}_q(n_y),
\label{lem:omegdef}
\end{eqnarray}
where $n_y, n_z$ are from
Definition
\ref{def:n}
and 
where $\Psi$ is from Definition
\ref{def:psi}.
We view $\Omega$ as a linear operator that acts on
finite dimensional 
$U_q(\mathfrak{sl}_2)$-modules.
\end{definition}

\noindent We now present our main result.
\begin{theorem}
\label{thm:main}
For the operator $\Omega $ from Definition
\ref{def:omega} the following hold
on each finite dimensional 
$U_q(\mathfrak{sl}_2)$-module:
\begin{eqnarray*}
\Omega^{-1}x\, \Omega = y,
\qquad \qquad 
\Omega^{-1}y\,  \Omega = z,
\qquad \qquad 
\Omega^{-1}z\, \Omega = x.
\end{eqnarray*}
\end{theorem}
\noindent {\it Proof:} 
Observe 
\begin{eqnarray*}
\Omega^{-1}x\, \Omega
&=&
{\rm{exp}}_q(n_y)^{-1}\, \Psi^{-1} \,
{\rm{exp}}_q(n_z)^{-1}\,x 
\,{\rm{exp}}_q(n_z)\, \Psi \,
{\rm{exp}}_q(n_y)
\\
&=& {\rm{exp}}_q(n_y)^{-1}\, \Psi^{-1} 
xyx 
 \Psi \,
{\rm{exp}}_q(n_y) \qquad  \qquad \qquad \qquad  \qquad 
(\mbox{by Lemma
\ref{lem:xyx}(i)})
\\
&=& {\rm{exp}}_q(n_y)^{-1}\, \Psi^{-1}
xy 
 \Psi  x\,
{\rm{exp}}_q(n_y) \qquad \qquad \qquad \qquad \qquad (\mbox{by Lemma 
\ref{lem:psiprop}(i)})
\\
&=& {\rm{exp}}_q(n_y)^{-1}\, \Psi^{-1}
(1-q^{-1}(q-q^{-1})n_z) 
 \Psi  x\,
{\rm{exp}}_q(n_y) \qquad \;
(\mbox{by 
(\ref{eq:nz})})
\\
&=& {\rm{exp}}_q(n_y)^{-1}\,(x-q^{-1}(q-q^{-1})x^{-1}n_z)
\,{\rm{exp}}_q(n_y) \qquad \qquad 
(\mbox{by Lemma
\ref{lem:psiprop}(iii)})
\\
&=& {\rm{exp}}_q(n_y)^{-1}\,(x-x^{-1}+y)
\,{\rm{exp}}_q(n_y) \qquad \qquad 
\qquad \qquad \;\;
(\mbox{by 
(\ref{eq:nz})})
\\
&=& y 
\qquad \qquad \qquad \qquad \qquad \qquad \qquad \qquad
\qquad \qquad \qquad \;\;
(\mbox{by Lemma 
\ref{lem:xyy2}(ii)})
\end{eqnarray*}
\noindent and 
\begin{eqnarray*}
\Omega^{-1}y\, \Omega
&=&
{\rm{exp}}_q(n_y)^{-1}\, \Psi^{-1} \,
{\rm{exp}}_q(n_z)^{-1}\,y
\,{\rm{exp}}_q(n_z)\, \Psi \,
{\rm{exp}}_q(n_y)
\\
&=& 
{\rm{exp}}_q(n_y)^{-1}\, \Psi^{-1} x^{-1}
\Psi \,
{\rm{exp}}_q(n_y)
\qquad \qquad \qquad \qquad 
(\mbox{by Lemma 
\ref{lem:expform2}(ii)})
\\
&=& 
{\rm{exp}}_q(n_y)^{-1}\, x^{-1}
\,
{\rm{exp}}_q(n_y)
\;\qquad \qquad \qquad \qquad  \qquad (\mbox{by Lemma
\ref{lem:psiprop}(i)})
\\
&=& 
z
\qquad \qquad \qquad \qquad \qquad 
\qquad  \qquad \qquad \qquad \qquad 
(\mbox{by Lemma
\ref{lem:expform2}(i)})
\end{eqnarray*}
\noindent and
\begin{eqnarray*}
\Omega^{-1}z\, \Omega
&=&
{\rm{exp}}_q(n_y)^{-1}\, \Psi^{-1} \,
{\rm{exp}}_q(n_z)^{-1}\,z
\,{\rm{exp}}_q(n_z)\, \Psi \,
{\rm{exp}}_q(n_y)
\\
&=& 
{\rm{exp}}_q(n_y)^{-1}\, \Psi^{-1} (z+x-x^{-1}) 
\Psi \,
{\rm{exp}}_q(n_y)
\qquad \qquad \qquad \quad 
(\mbox{by Lemma
\ref{lem:xyy}(iii)})
\\
&=& 
{\rm{exp}}_q(n_y)^{-1}\, \Psi^{-1} (x-q^{-1}(q-q^{-1})n_yx^{-1}) 
\Psi \,
{\rm{exp}}_q(n_y)
\qquad   \; (\mbox{by 
(\ref{eq:ny})})
\\
&=& 
{\rm{exp}}_q(n_y)^{-1}\,  (x-q^{-1}(q-q^{-1}) \Psi^{-1} n_y \Psi x^{-1}) 
\,
{\rm{exp}}_q(n_y)
\qquad   \; 
(\mbox{by Lemma
\ref{lem:psiprop}(i)})
\\
&=& 
{\rm{exp}}_q(n_y)^{-1}\, (x-q^{-1}(q-q^{-1})x n_y) 
\,{\rm{exp}}_q(n_y)
\qquad \qquad \qquad  (\mbox{by Lemma
\ref{lem:psiprop}(ii)})
\\
&=& 
{\rm{exp}}_q(n_y)^{-1}\, xzx
\,{\rm{exp}}_q(n_y)
\; \;\qquad \qquad \qquad \qquad \qquad \qquad \quad (\mbox{by 
(\ref{eq:ny})})
\\
&=& 
x
\qquad \qquad 
\qquad \qquad \qquad \qquad \qquad \qquad  
 \qquad \quad \;\qquad \qquad  
(\mbox{by Lemma 
\ref{lem:xyx2}(i)}).
\end{eqnarray*}
\hfill $\Box $ \\

\noindent We finish this section with a comment.

\begin{corollary}
On a finite dimensional 
 $U_q(\mathfrak{sl}_2)$-module, $\Omega^3$
 commutes with the action of each element of
 $U_q(\mathfrak{sl}_2)$.
\end{corollary}
\noindent {\it Proof:} 
Immediate from
Theorem
\ref{thm:main} and since 
$x^{\pm 1}, y,z$ generate
 $U_q(\mathfrak{sl}_2)$.
\hfill $\Box $ \\

\section{The action of $\Omega$ on $L(n,\varepsilon)$}

\noindent In this section we describe the action of $\Omega$ on the module
$L(n,\varepsilon)$. We will do this by displaying the action
of $\Omega$ and $\Omega^{-1}$ on the basis for
$L(n,\varepsilon)$ given in
Lemma \ref{lem:altact}. 
We begin with a few observations.

\begin{lemma}
\label{lem:psiln}
For an integer $n\geq 0$ and for $\varepsilon \in \lbrace 1, -1\rbrace$
  let $u_0, u_1, \ldots, 
u_n$ denote the basis for $L(n,\varepsilon)$ given 
in Lemma
\ref{lem:altact}. Then 
$\Psi u_i = q^{2i(n-i)+(s-n^2)/2} u_i$ for $0 \leq i \leq n$,
where $s=0$ (if $n$ is even) and $s=1$ (if $n$ is odd).
\end{lemma}
\noindent {\it Proof:} Immediate
from Definition
\ref{def:psi}.
\hfill $\Box $ \\

\begin{lemma}
\label{lem:nynzln}
For an integer $n\geq 0$ and for $\varepsilon \in \lbrace 1, -1\rbrace$
  let $u_0, u_1, \ldots, 
u_n$ denote the basis for $L(n,\varepsilon)$ given 
in Lemma
\ref{lem:altact}. Then
the following 
{\rm (i), (ii)} hold.
\begin{enumerate}
\item[{\rm (i)}] 
$n_y u_i = -q^{n-i}\lbrack n-i+1 \rbrack u_{i-1}
\;(1 \leq i \leq n),  \quad n_y u_0=0$.
\item[{\rm (ii)}] 
$n_z u_i = q^{-i} \lbrack i+1 \rbrack u_{i+1} 
\;(0 \leq i \leq n-1), \quad n_z u_n=0$.
\end{enumerate}
\end{lemma}
\noindent {\it Proof:} 
Use 
Lemma \ref{lem:altact}
and
Definition
\ref{def:n}.
\hfill $\Box $ \\

\noindent 
We recall some notation.
For  integers $n\geq i\geq 0$
we define
\begin{eqnarray*}
\biggl\lbrack {{n}\atop{i}} \biggr\rbrack
= \frac{\lbrack n \rbrack^!}{\lbrack i \rbrack^! \lbrack n-i \rbrack^!}.
\end{eqnarray*}

\begin{lemma}
\label{lem:expn}
For an integer $n\geq 0$ and for $\varepsilon \in \lbrace 1, -1\rbrace$
  let $u_0, u_1, \ldots, 
u_n$ denote the basis for $L(n,\varepsilon)$ given 
in Lemma
\ref{lem:altact}. Then
for $0 \leq j \leq n$ we have
\begin{eqnarray}
{\rm{exp}}_q  (n_y)\,u_j &=& \sum_{i=0}^j (-1)^{i+j} q^{(j-i)(n-i-1)} 
\biggl\lbrack {{n-i}\atop{j-i}} \biggr\rbrack
u_i,
\label{eq:expny}
\\
{\rm{exp}}_q(n_y)^{-1}\,u_j &=& \sum_{i=0}^j q^{(j-i)(n-j)} 
\biggl\lbrack {{n-i}\atop{j-i}} \biggr\rbrack
u_i,
\label{eq:expnyinv}
\\
{\rm{exp}}_q(n_z)\,u_j &=& \sum_{i=j}^{n} q^{j(j-i)} 
\biggl\lbrack {{i}\atop{j}} \biggr\rbrack
u_i,
\label{eq:expnz}
\\
{\rm{exp}}_q(n_z)^{-1}\,u_j &=& \sum_{i=j}^{n} (-1)^{i+j} q^{(j-i)(i-1)} 
\biggl\lbrack {{i}\atop{j}} \biggr\rbrack
u_i.
\label{eq:expnzinv}
\end{eqnarray}
\end{lemma}
\noindent {\it Proof:}
In order to verify 
(\ref{eq:expny}) and 
(\ref{eq:expnz}),
evaluate
the left-hand side
using
Lemma \ref{lem:nynzln}
and
Definition 
\ref{def:tnil}.
Lines 
(\ref{eq:expnyinv}) and 
(\ref{eq:expnzinv})
are similarly verified using
Lemma \ref{lem:expinv}.
\hfill $\Box $ \\

\begin{theorem}
\label{thm:omegamat}
For an integer $n \geq 0$ and for
$\varepsilon \in \lbrace 1,-1\rbrace $ let $u_0, u_1, \ldots, 
u_n$ denote the basis for $L(n,\varepsilon)$ given 
in Lemma
\ref{lem:altact}. Then for $0 \leq j\leq n$ we have
\begin{eqnarray}
\Omega u_j &=& 
\sum_{i=0}^{n-j} (-1)^j q^{(n-i-1)j+(s-n^2)/2} 
\biggl\lbrack {{n-i}\atop{j}} \biggr\rbrack
u_i,
\label{eq:omeguj}
\\
\Omega^{-1} u_j &=& 
\sum_{i=n-j}^n (-1)^{n-j} q^{(1-i)(n-j)+(n^2-s)/2} 
\biggl\lbrack {{i}\atop{n-j}} \biggr\rbrack
u_i,
\label{eq:omegujinv}
\end{eqnarray}
where $s=0$ (if $n$ is even) and $s=1$ (if $n$ is odd).
\end{theorem}
\noindent {\it Proof:}
In order to verify 
(\ref{eq:omeguj}),
evaluate the left-hand side using
(\ref{lem:omegdef}),
Lemma
\ref{lem:psiln},
(\ref{eq:expny}), 
(\ref{eq:expnz}), and simplify the result using
the $q$-Vandermonde summation formula
\cite[p.~11]{GR}.
Line
(\ref{eq:omegujinv}) is similarly verified.
\hfill $\Box $ \\

\noindent
We finish this section with a comment.

\begin{corollary}
For an integer $n \geq 0$ and for
$\varepsilon \in \lbrace 1,-1\rbrace $,
$\Omega^3$ acts as a scalar multiple of the identity on
$L(n,\varepsilon)$. The scalar is
$q^{-n(n+2)/2}$ (if $n$ is even) and
$-q^{(1-n)(n+3)/2}$ (if $n$ is odd).
\end{corollary}
\noindent {\it Proof:} 
Routine calculation using
Theorem \ref{thm:omegamat}.
\hfill $\Box $ \\


\section{Remarks}

\medskip
\noindent In this section we make some remarks and tie
up some loose ends.

\begin{remark}
\label{rem:Fairlie}
\rm
In \cite{fairlie} Fairlie considers an
associative $\K$-algebra
with generators $X,Y,Z$ and relations
\begin{eqnarray}
qXY-q^{-1}YX &=& Z,
\label{eq:cyc1}
\\
qYZ-q^{-1}ZY &=& X,
\label{eq:cyc2}
\\
qZX-q^{-1}XZ &=& Y.
\label{eq:cyc3}
\end{eqnarray}
He interprets this algebra as a 
$q$-deformation of $\mbox{SU}(2)$ and 
he works out the irreducible representations.
See
\cite[Remark 8.11]{ciccoli},
\cite[Section 3]{flat},
\cite{havlicek},
\cite{odesskii} for related work.
In spite of the superficial resemblance we do not see
any connection between
(\ref{eq:cyc1})--(\ref{eq:cyc3})
and the equitable presentation
of 
$U_q(\mathfrak{sl}_2)$.
\end{remark}

\begin{remark}
\rm 
In 
\cite{Zhidd}
A.~S.~Zhedanov
introduced the Askey-Wilson algebra.
He 
used it to study the Askey-Wilson polynomials
and related polynomials in the Askey scheme
\cite{KoeSwa}.
The following attractive
version of the algebra appears in
\cite[p.~101]{rosen},
\cite{smithandbell},
\cite[Section 3.3.3]{wieg}.
For a sequence of scalars
$g_x, g_y, g_z, h_x, h_y, h_z$ taken from
$\K$, the corresponding Askey-Wilson algebra
is the unital associative $\K$-algebra
with generators $X,Y,Z$ and relations
\begin{eqnarray}
\label{eq:wieg1}
qXY-q^{-1}YX &=& g_zZ+h_z,
\\
\label{eq:wieg2}
qYZ-q^{-1}ZY &=& g_xX+h_x,
\\
\label{eq:wieg3}
qZX-q^{-1}XZ &=& g_yY+h_y.
\end{eqnarray}
See
\cite{GYZnature},
\cite{GYLZmut},
\cite{GYZTwisted},
  \cite{GYZlinear},
\cite{GYZspherical},
\cite{aw},
\cite{ZheCart},
\cite{Zhidden}
for work involving
the Askey-Wilson algebra.  
We note that for $g_x=g_y=g_z=1$ and $h_x=h_y=h_z=0$
the relations
(\ref{eq:wieg1})--(\ref{eq:wieg3}) become
(\ref{eq:cyc1})--(\ref{eq:cyc3}).
Moreover for $g_x=g_y=g_z=0$ and $h_x=h_y=h_z=q-q^{-1}$
the relations
(\ref{eq:wieg1})--(\ref{eq:wieg3}) become
(\ref{eq:2buq2})--(\ref{eq:2buq4}).
In this case, and referring to 
Theorem \ref{thm:uq2},
the Askey-Wilson algebra is isomorphic
to the subalgebra of 
$U_q(\mathfrak{sl}_2)$ generated by
$x,y,z$. As far as we know, this case of
the Askey-Wilson algebra has not yet been considered
by other researchers.
\end{remark}

\begin{remark}
\rm
In the literature one can find many
presentations of algebras 
that are related in some way to 
$\mathfrak{sl}_2$.
See for example
\cite{bavula},
\cite{bavoy},
\cite{du},
\cite{cas},
\cite{delb},
\cite{fuj},
\cite{bruyn},
\cite{bruyn2},
\cite{woro},
\cite[p. 48]{rosen},
\cite{spsmith},
\cite{van},
\cite{witten1},
\cite{witten}.
As far as we know,
none of these has a direct connection to
the equitable presentation of
$U_q(\mathfrak{sl}_2)$.
\end{remark}

\begin{remark}
\rm
The equitable presentation for
$U_q(\mathfrak{sl}_2)$ 
appears  implicitly as part of  
a presentation given in
\cite[Theorem 2.1]{tdanduq}
for the quantum affine
algebra
$U_q(\widehat{\mathfrak{sl}}_2)$.
\end{remark}

\begin{remark}
\rm In light of
Theorem
\ref{thm:uq2} and Corollary
\ref{cor:yzinv} it is natural to consider a unital
associative $\K$-algebra 
that has generators
$x^{\pm 1}$, $y^{\pm 1}$, $z^{\pm 1}$
and relations
\begin{eqnarray*}
&&xx^{-1} = 
x^{-1}x =  1,
\qquad \qquad 
yy^{-1} = 
y^{-1}y =  1,
\qquad \qquad 
zz^{-1} = 
z^{-1}z =  1,
\qquad \qquad 
\\
&&\frac{qxy-q^{-1}yx}{q-q^{-1}}=1,
\qquad \qquad 
\frac{qyz-q^{-1}zy}{q-q^{-1}}=1,
\qquad \qquad 
\frac{qzx-q^{-1}xz}{q-q^{-1}}=1.
\end{eqnarray*}
We denote this algebra by  
$U^{\Delta}_q(\mathfrak{sl}_2)$ and call it the {\it equitable}
$q$-deformation of 
$\mathfrak{sl}_2$. 
We invite the reader to investigate
$U^{\Delta}_q(\mathfrak{sl}_2)$.
\end{remark}

\begin{remark}
\rm
For a symmetrizable Kac-Moody  Lie  algebra
$\mathfrak{g}$ the second author
has obtained 
a presentation for  
the quantum group $U_q(\mathfrak{g})$
that is analogous to 
the equitable presentation 
for $U_q(\mathfrak{sl}_2)$.
See \cite{km} for the details. 
\end{remark}

\begin{remark}
\rm
We discovered the equitable presentation for
$U_q(\mathfrak{sl}_2)$ 
during our recent study of tridiagonal pairs
\cite{TD00},
\cite{shape},
\cite{tdanduq},
and the closely related Leonard pairs
\cite{LS99},
  \cite{qSerre},
   \cite{LS24},
   \cite{conform},
    \cite{lsint},
\cite{TLT:split},
 \cite{TLT:array},
\cite{qrac},
\cite{aw}.
A Leonard pair is a pair of semi-simple
linear transformations on a finite-dimensional vector
space, each of which acts tridiagonally on an eigenbasis
for the other 
\cite[Definition 1.1]{LS99}. There is a close connection between
Leonard pairs and the orthogonal polynomials that make
up the terminating branch of the Askey scheme 
\cite{KoeSwa},
\cite[Appendix A]{LS99},
 \cite{TLT:array}.
 A tridiagonal pair is a mild generalization of
a Leonard pair
\cite[Definition 1.1]{TD00}.
\end{remark}


\noindent Tatsuro Ito \hfil\break
\noindent Department of Computational Science \hfil\break
\noindent Faculty of Science \hfil\break
\noindent Kanazawa University \hfil\break
\noindent Kakuma-machi \hfil\break
\noindent Kanazawa 920-1192, Japan \hfil\break
\noindent Email:  {\tt ito@kappa.s.kanazawa-u.ac.jp}

\bigskip   

\noindent Paul Terwilliger \hfil\break
Department of Mathematics \hfil\break
University of Wisconsin \hfil\break
480 Lincoln Drive \hfil\break
Madison, Wisconsin, 53706 USA 
\hfil\break
{Email: \tt terwilli@math.wisc.edu} 

\bigskip

\noindent Chih-wen Weng \hfil\break
Department of Applied Mathematics \hfil\break
National Chiao Tung University  \hfil\break
1001 Ta Hsueh Road \hfil\break
Hsinchu 30050, Taiwan, ROC 
\hfil\break
{Email: \tt weng@math.nctu.edu.tw} \hfil\break

\end{document}